\documentclass[12pt]{article}

\usepackage{amssymb,amsmath,amsthm, amsfonts}

\def\sqr#1#2{{\vcenter{\vbox{\hrule height.#2pt
              \hbox{\vrule width.#2pt height#1pt \kern#1pt \vrule width.#2pt}
              \hrule height.#2pt}}}}
\def\signed #1{{\unskip\nobreak\hfil\penalty50
              \hskip2em\hbox{}\nobreak\hfil#1
              \parfillskip=0pt \finalhyphendemerits=0 \par}}
\def\endpf{\signed {$\sqr69$}}

\def\dbR{{\mathop{\rm l\negthinspace R}}}

\def\3n{\negthinspace \negthinspace \negthinspace }
\def\2n{\negthinspace \negthinspace }
\def\1n{\negthinspace }

\def\see{{\it see} }
\def\ds{\displaystyle}

\def\dbN{{\mathop{\rm l\negthinspace N}}}

\def\dbR{{\mathop{\rm l\negthinspace R}}}
\def\={\buildrel \triangle \over =}

\def\resp{{\it resp. }}
\def\a{\alpha}

\def\l{\lambda}

 \def\n{\nabla}

\def\t{\times}
\def\f{\varphi}
\def\th{\theta}
\def\o{\omega}

\def\i{\infty}
\def\ns{\noalign{\ss} }

\def\O{\Omega}

\def\cF{{\cal F}}

\def\cO{{\cal O}}

\def\cl{{\cal l}}
\def\no{\noindent}

\def\ss{\smallskip}
\def\ms{\medskip}
\def\bs{\bigskip}
\def\q{\quad}
\def\qq{\qquad}
\def\hb{\hbox}

\def\lan{\mathop{\langle}}
\def\ran{\mathop{\rangle}}

\def\min{\mathop{\rm min}}
\def\exp{\mathop{\rm exp}}

\def\pa{\partial}

\def\cd{\cdot}

\def\as{\hbox{\rm a.s.{ }}}

\def\supp{\hbox{\rm supp$\,$}}

\def\cl{\overline}

\def\|{\Big |}
\def\({\Big (}
\def\){\Big )}
\def\[{\Big[}
\def\]{\Big]}
\def\be{\begin{equation}}
\def\bel{\begin{equation}\label}
\def\ee{\end{equation}}
\def\bt{\begin{theorem}}
\def\bcd{\begin{condition}}
\def\ecd{\end{condition}}
\def\et{\end{theorem}}
\def\bc{\begin{corollary}}
\def\ec{\end{corollary}}
\def\bde{\begin{definition}}
\def\ede{\end{definition}}
\def\bl{\begin{lemma}}
\def\el{\end{lemma}}
\def\bp{\begin{proposition}}
\def\ep{\end{proposition}}
\def\br{\begin{remark}}
\def\er{\end{remark}}
\def\ba{\begin{array}}
\def\ea{\end{array}}
\def\ed{\end{document}}
\def\ns{\noalign{\ms}}
\def\ds{\displaystyle}

\def\square#1{\vbox{\hrule\hbox{\vrule height#1%
     \kern#1\vrule}\hrule}}
\def\rectangle#1#2{\vbox{\hrule\hbox{\vrule height#1%
     \kern#2\vrule}\hrule}}

\font\tenbb=msbm10 \font\sevenbb=msbm7 \font\fivebb=msbm5

\newfam\bbfam
\scriptscriptfont\bbfam=\fivebb \textfont\bbfam=\tenbb
\scriptfont\bbfam=\sevenbb

\newtheorem{lemma}{Lemma}
\newtheorem{remark}{Remark}

\newtheorem{theorem}{Theorem}

\begin{document}
\title{\bf Unique Continuation for Stochastic Parabolic Equations}

\author{Xu Zhang\thanks{Academy of Mathematics and Systems Sciences, Academia Sinica,
Beijing 100080, China; and Yangtze Center of Mathematics, Sichuan
University, Chengdu 610064, China. {\small\it E-mail:} {\small\tt
xuzhang@amss.ac.cn}. This work was supported by the NSF of China
under grants 10371084 and 10525105, the NCET of China under grant
NCET-04-0882. }}

\date{}

\maketitle

\begin{abstract}
\no This paper is devoted to a study of the unique continuation
property for stochastic parabolic equations. Due to the adapted
nature of solutions in the stochastic situation, classical
approaches to treat the the unique continuation problem for
deterministic equations do not work. Our method is based on a
suitable partial Holmgren coordinate transform and a stochastic
version of Carleman-type estimate.

\end{abstract}

\bs

\no{\bf 2000 Mathematics Subject Classification}. Primary: 60H15;
Secondary: 34A12.

\bs

\no{\bf Key Words}.  Unique Continuation,  stochastic parabolic
equations, Carleman-type estimate, partial Holmgren coordinate
transform.

\newpage

\section{Introduction and main result}

Let $T>0$, $G \subset \dbR^{n}$ ($n\in\dbN$) be a given bounded
domain with a $C^2$ boundary $\pa G$, and $G_0\neq G$ be a given
subdomain of $G$. Put $
 Q\=(0,T)\times G$ and $Q_0\=(0,T)\times G_0$.
Throughout this paper, we assume that $a^{ij}\in
W^{1,\i}(0,T;W_{loc}^{2,\i}(G))$ satisfy $a^{ij}=a^{ji}$
($i,j=1,2,\cdots,n$) and for any open subset $G_1$ of $G$, there is
a constant $s_0=s_0(G_1)>0$ so that
 \bel{h1}
 \sum_{i,j}a^{ij}\xi^{i}\xi^{j}
 \geq s_0|\xi|^{2},\qq \forall\;(t,x,\xi)\equiv (t,x,\xi^{1},\cdots,\xi^{n})
 \in(0,T)\times G_1\t \dbR^{n}.
 \ee
Here, we denote $\ds\sum^{n}_{i,j=1}$ simply by $\ds\sum_{i,j}$.

Let $(\O,\cF,\{\cF_t\}_{t\ge 0},P)$  be a complete filtered
probability space on which a $1$ dimensional standard Brownian
motion $\{w(t)\}_{t\ge0}$ is defined. Let $H$ be a Fr\'echet space.
We denote by $L_{\cF}^2(0,T;H)$ the Fr\'echet space consisting of
all $H$-valued $\{\cF_t\}_{t\ge 0}$-adapted processes $X(\cd)$ such
that $\mathbb{E}(|X(\cd)|_{L^2(0,T;H)}^2)<\i$, with the canonical
quasi-norm; by $L_{\cF}^\i(0,T;H)$ the Fr\'echet space consisting of
all $H$-valued $\{\cF_t\}_{t\ge 0}$-adapted bounded processes, with
the canonical quasi-norm; and by $L_{\cF}^2(\O;C([0,T];H))$ the
Fr\'echet space consisting of all $H$-valued $\{\cF_t\}_{t\ge
0}$-adapted continuous processes $X(\cd)$ such that
$\mathbb{E}(|X(\cd)|_{C([0,T];H)}^2)<\i$, with the canonical
quasi-norm.

Let us consider the following stochastic parabolic equation:
 \bel{hh6.1}
 \cF z\equiv dz-\sum_{i,j}(a^{ij}z_i)_jdt=[\lan a,\n z\ran+bz]dt+cz
 dw(t)\qq\hb{ in }Q.
 \ee
Here $a,b$ and $c$ are  suitable coefficients. For simplicity, we
use the notation $ z_i\equiv z_i(x)={{\pa z(x)}/{\pa x_i}}, $ where
$x_i$ is the $i$-th coordinate of a generic point
$x=(x_1,\cdots,x_n)$ in $\dbR^n$. In a similar manner, in the sequel
we use the notation $ u_i$, $v_i$, etc. for the partial derivatives
of $u$ and $v$ with respect to $x_i$. Also, we denote the scalar
product in $\dbR^n$ by $\lan\cd,\cd\ran$.

The main result of this paper is stated as follows:

\bt\label{1t3}
Let $a\in L^{\i}_{\cF}(0,T;L_{loc}^{\i}(G;\dbR^n))$, $b\in
L^{\i}_{\cF}(0,T;L_{loc}^{\i}(G))$, and $c\in L^{\i}_{\cF}(0,T;$
$W_{loc}^{1,\i}(G))$. Then any solution $z\in
L_{\cF}^2(\O;C([0,T];L_{loc}^2(G)))\bigcap
L_{\cF}^2(0,T;H_{loc}^1(G))$  of (\ref{hh6.1}) vanishes identically
in $Q\t\O, \as dP$ provided that $z=0$ in $Q_0\t\O, \as dP$.
\et

The above result is a unique continuation theorem for stochastic
parabolic equations. There are numerous references on the unique
continuation for deterministic parabolic equations (see, for
example, \cite{I, SS, S, Y} and so on). However, to the author's
best acknowledge, nothing is known for its stochastic counterpart.

There are two classical tools in the study of the unique
continuation for deterministic partial differential equations. One
is Holmgren-type uniqueness theorem, another is  Carleman-type
estimate. Note however that the solution of a stochastic equation is
generally non-analytic in time even if the coefficients of the
equation are constants. Therefore, one cannot expect a Holmgren-type
uniqueness theorem for the unique continuation for stochastic
equations except some very special cases. On the other hand, the
usual approach to employ Carleman-type estimate for the unique
continuation needs to localize the problem. The difficulty of our
present stochastic problem consists in the fact that one cannot
simply localize the problem as usual because the classical
localization technique may change the adaptedness of solutions,
which is a key feature in the stochastic setting. In our equation
(\ref{hh6.1}), for the space variable $x$, we may proceed as in the
classical argument. However, for the time variable $t$, due to the
adaptedness requirement, we will have to treat it in a deliberate
way. For this purpose, we shall introduce a suitable ``partial
Holmgren coordinate transform" (\see (\ref{pa1})) and deduce a key
stochastic version of Carleman-type estimate (\see Theorem
\ref{c1t3} in the next section).

It is well-known that, unique continuation is an important problem
not only in partial differential equations itself, but also in some
application problems such as controllability (\cite{Zu}), inverse
problems (\cite{IS}), optimal control (\cite{LY}) and so on.
Numerous studies on unique continuation for deterministic partial
differential equations can be found in \cite{Ho2,Zui} and the rich
references cited therein. It would be quite interesting to extend
the deterministic unique continuation theorems to the stochastic
ones, but there are many things which remain to be done, and some of
which seem to be challenging. In this paper, in order to present the
key idea in the simplest way, we do not pursue the full technical
generality.

The rest of this paper is organized as follows. In Section 2, as a
key preliminary, we show a Carleman-type estimate for stochastic
parabolic operators. Section 3 is devoted to the proof of Theorem
\ref{1t3}.

\section{Carleman-type estimate for stochastic  parabolic operators}\label{s5}

For any nonnegative and nonzero function $\psi\in C^3(\cl{G})$,
any $k\ge 2$, and any (large) parameters $\l>1$ and $\mu>1$, put
 \be
 \label{alphad}
 \ell=\l\a,\q\alpha(t,x)={e^{\mu\psi(x)}-e^{2\mu|\psi|_{C(\cl{G })}}\over t^k(T-t)^k},\q
 \varphi(t,x)={e^{\mu\psi(x)}\over t^k(T-t)^k}.
 \ee

In the sequel, we will use $C$ to denote a generic positive
constant depending only on $T$, $G$, $G_0$ and $(a^{ij})_{n\t n}$,
which may change from line to line. Also, for $r\in \dbN$, we
denote by $O(\mu^r)$ a function of order $\mu^r$ for large $\mu$
(which is independent of $\l$); by $O_{\mu}(\l^r)$ a function of
order $\l^r$ for fixed $\mu$ and for large $\l$. We recall the
following known result.

\bl\label{c1t2}
{\rm (\cite{TZ0, TZ})} Let $b^{ij}\in C^{1,2}(\cl{Q})$ satisfying
$b^{ij}=b^{ji}$. Assume that either  $(b^{ij})_{n\t n}$ or
$-(b^{ij})_{n\t n}$ is a uniformly positive definite matrix, and
$s_0(>0)$ is its smallest eigenvalue. Let $u$ be a
$C^2(\cl{G})$-valued semimartingale. Set
 \bel{h5}
 \th=e^{\ell }, \q v=\th u, \q \Psi=2\sum_{i,j}b^{ij}\ell_{ij}.
 \ee
Then for any $x\in G$ and $\o\in\O$ $(\as dP)$,
 \bel{c1e14}
 \ba{ll}
 \displaystyle
 2\int_0^T\th\[-\sum_{i,j} (b^{ij}v_i)_j+Av\]\[du-\sum_{i,j}(b^{ij}u_i)_jdt\]+2\int_0^T\sum_{i,j} (b^{ij}v_idv)_j\\
 \noalign{\ss}
 \displaystyle\q+ 2\int_0^T\sum_{i,j}\[\sum_{i',j'}\(2b^{ij} b^{i'j'}\ell_{i'}v_iv_{j'}
 -b^{ij}b^{i'j'}\ell_iv_{i'}v_{j'}\)+\Psi
 b^{ij}v_iv\\
 \noalign{\ss}
 \displaystyle\qq\qq\qq\qq\qq\qq- b^{ij}\(A\ell_i+{\Psi_i\over2}\)v^2\]_jdt\\
 \noalign{\ss}
 \displaystyle
 \ge 2\sum_{i,j}\int_0^Tc^{ij}v_iv_jdt
 +\int_0^TBv^2dt+\int_0^T\|-\sum_{i,j} (b^{ij}v_i)_j+Av\|^2dt\\
 \noalign{\ss}
 \displaystyle\q-\int_0^T\th^2\sum_{i,j}
 b^{ij}du_idu_j-\int_0^T\th^2\[A-\sum_{i,j}
 \(b^{ij}\ell_i\ell_j+(b^{ij}\ell_i)_j\)\](du)^2,
 \ea
 \ee
where
\bel{c1e15}
\3n\left\{
 \ba{ll}
 \ds A\=-\sum_{i,j}
 \[b^{ij}\ell_i\ell_j-(b^{ij}\ell_i)_j\]-\Psi,\\
  \ns
 \ds
 B\=2\[A\Psi-
 \sum_{i,j}(Ab^{ij}\ell_i)_j\] -A_t-\sum_{i,j} (b^{ij}\Psi_j)_i-\ell_t^2,\\
 \ns
 \ds c^{ij}\=\sum_{i',j'}\[2b^{ij'}(b^{i'j}\ell_{i'})_{j'} -
 (b^{ij}b^{i'j'}\ell_{i'})_{j'}\]-{b_t^{ij}\over2}+\Psi b^{ij}.
  \ea
\right.
 \ee
Moreover, for $\l$ and $\mu$ large enough, it holds
 \bel{h6}
 \ba{ll}\ds
 A=-\l^2\mu^2\f^2\sum_{i,j}b^{ij}\psi_i\psi_j
  +\l\f O(\mu^2), \\
 \ns
 \ds B\ge 2s_0^2\l^3\mu^4\f ^3|\n \psi|^4+\l^3\f^3O(\mu^3)
  +\l^2\f^2O(\mu^4)+\l\f O(\mu^4)\\ \ns
  \ds\qq +\l^2\f^{2+2k^{-1}}O(e^{4\mu
|\psi|_{C(\cl{G})}})+\l^2\f^{2+k^{-1}} O(\mu^2)+\l\f^{1+k^{-1}}
O(\mu^2), \\
 \ns
 \ds
 \sum_{i,j} c^{ij}v_iv_j\ge [s_0^2\l\mu^2\f |\n\psi|^2+\l\f  O(\mu)]|\n
 v|^2.
  \ea
 \ee
\el

We now show a Carleman-type estimate  for stochastic parabolic
operators as follows:

\bt\label{c1t3}
Let $b^{ij}\in C^{1,2}(\cl{Q})$ satisfying $b^{ij}=b^{ji}$. Assume
that either  $(b^{ij})_{n\t n}$ or $-(b^{ij})_{n\t n}$ is a
uniformly positive definite matrix. Let $\psi\in C^3(\cl{G})$
satisfy
 \bel{u1}
 \min_{x\in G}|\n\psi(x)|>0.
 \ee
Then there is some $\mu_0>0$ such that for all $\mu\ge \mu_0$, one
can find two constants $C=C(\mu)>0$ and $\l_1=\l_1(\mu)$ so that
for all $u\in L_{\cF}^2(\O;C([0,T];L^2(G)))\bigcap $ $
L_{\cF}^2(0,T;H_0^2(G))$, $f\in L_{\cF}^2(0,T;L^2(G))$ and $g\in
L_{\cF}^2(0,T; H^1(G))$ with
 \bel{h9}
 du-\sum_{i,j}(b^{ij}u_i)_jdt=fdt+gdw(t), \qq\hb{ in }Q,
 \ee
and all $\l\ge \l_1$, it holds
 \bel{h7}
 \ba{ll}\ds
 \l^3\mu^4\mathbb{E}\int_Q\f^3\th^2u^2dxdt+\l\mu^2\mathbb{E}\int_Q\f\th^2|\n
 u|^2dxdt\\
  \ns
  \ds
 \displaystyle\le C\Big\{\mathbb{E}\int_Q \th^2f^2dxdt+
 \mathbb{E}\int_Q \th^2\sum_{i,j}
 b^{ij}g_ig_jdxdt\\
 \ns
 \ds\q+\mathbb{E}\int_Q \th^2\[A-\sum_{i,j}
 \(b^{ij}\ell_i\ell_j+(b^{ij}\ell_i)_j\)\] g^2dxdt\Big\},
 \ea
 \ee
where
 \bel{h23}
 A\=-\sum_{i,j}
 (b^{ij}\ell_i\ell_j-b^{ij}_j\ell_i
 +b^{ij}\ell_{ij}).
 \ee
\et

{\it Proof.} Recalling that $k\ge 2$ and (\ref{alphad}), we get
 \bel{z01}
  \ba{ll}\ds
 |\l^2\f^2O(\mu^4)+\l\f O(\mu^4)+\l^2\f^{2+2k^{-1}}O(e^{4\mu
|\psi|_{C(\cl{G})}})+\l^2\f^{2+k^{-1}} O(\mu^2)+\l\f^{1+k^{-1}}
O(\mu^2)|\\
 \ns
\le \f^3O_{\mu}(\l^2).
 \ea
 \ee
Integrating (\ref{c1e14}) (in Lemma \ref{c1t2}) on $G$, taking
mean value in both sides, and noting (\ref{h6}) (in Lemma
\ref{c1t2}) and (\ref{z01}), recalling that $u$, and hence $v$,
belongs to $ L_{\cF}^2(0,T;H_0^2(G))$, we conclude that there is a
constant $c_0>0$ such that
 \bel{h10}
 \ba{ll}
 \displaystyle
 2\mathbb{E}\int_Q\th\[-\sum_{i,j} (b^{ij}v_i)_j+Av\]\[du-\sum_{i,j}(b^{ij}u_i)_jdt\]dx\\
 \noalign{\ss}
 \displaystyle
 \ge 2c_0\mathbb{E}\int_Q[\l\mu^2\f |\n\psi|^2+\l\f  O(\mu)]|\n
 v|^2dtdx\\
 \noalign{\ss}
 \displaystyle\q
 +2c_0\mathbb{E}\int_Q\[\l^3\mu^4\f ^3|\n \psi|^4+\l^3\f^3O(\mu^3)
  +\f^3O_\mu(\l^2)\]v^2dtdx\\
 \noalign{\ss}
 \displaystyle\q+\mathbb{E}\int_Q\|-\sum_{i,j} (b^{ij}v_i)_j+Av\|^2dtdx\\
 \ns
 \ds\q-\mathbb{E}\int_Q\th^2\sum_{i,j}
 b^{ij}du_idu_jdx-\mathbb{E}\int_Q\th^2\[A-\sum_{i,j}
 \(b^{ij}\ell_i\ell_j+(b^{ij}\ell_i)_j\)\](du)^2dx.
 \ea
 \ee

By (\ref{h9}), we have
 \bel{h11}
 \3n\ba{ll}
 \displaystyle
 2\mathbb{E}\int_Q\th\[-\sum_{i,j} (b^{ij}v_i)_j+Av\]\[du-\sum_{i,j}(b^{ij}u_i)_jdt\]dx\\
 \noalign{\ss}
 \displaystyle=2\mathbb{E}\int_Q\th\[-\sum_{i,j} (b^{ij}v_i)_j+Av\][fdt+gdw(t)]dx\\
 \noalign{\ss}
 \displaystyle
 =2\mathbb{E}\int_Q\th\[-\sum_{i,j}
 (b^{ij}v_i)_j+Av\]fdtdx\\
 \noalign{\ss}
 \displaystyle\le \mathbb{E}\int_Q\|-\sum_{i,j}
 (b^{ij}v_i)_j+Av\|^2dtdx+\mathbb{E}\int_Q\th^2 f^2dtdx,
 \ea
 \ee
and
 \bel{h15}
 \ba{ll}\ds
 \mathbb{E}\int_Q\th^2\sum_{i,j}
 b^{ij}du_idu_jdx+\mathbb{E}\int_Q\th^2\[A-\sum_{i,j}
 \(b^{ij}\ell_i\ell_j+(b^{ij}\ell_i)_j\)\](du)^2dx\\
  \ns
 \ds= \mathbb{E}\int_Q \th^2\sum_{i,j}
 b^{ij}g_ig_jdxdt+\mathbb{E}\int_Q \th^2\[A-\sum_{i,j}
 \(b^{ij}\ell_i\ell_j+(b^{ij}\ell_i)_j\)\] g^2dxdt.
  \ea
 \ee

Combining (\ref{h10})--(\ref{h15}), we arrive at
 \bel{h16}
 \ba{ll}
 \displaystyle 2c_0\mathbb{E}\int_Q\f [\l\mu^2|\n\psi|^2+\l O(\mu)]|\n
 v|^2dtdx\\
 \noalign{\ss}
 \displaystyle\q
 +2c_0\mathbb{E}\int_Q\f^3\[\l^3\mu^4|\n \psi|^4+\l^3O(\mu^3)
  +O_\mu(\l^2)\]v^2dtdx\\
 \noalign{\ss}
 \displaystyle\le
 \mathbb{E}\int_Q\th^2f^2dxdt+
 \mathbb{E}\int_Q \th^2\sum_{i,j}
 b^{ij}g_ig_jdxdt\\
 \ns
 \ds\q+\mathbb{E}\int_Q \th^2\[A-\sum_{i,j}
 \(b^{ij}\ell_i\ell_j+(b^{ij}\ell_i)_j\)\] g^2dxdt.
 \ea
 \ee
Finally, combining  (\ref{h16}) and (\ref{u1}), and returning $v$
to $u$, we obtain the desired estimate (\ref{h7}).\endpf

\section{Proof of Theorem \ref{1t3}}\label{s6}

The proof is divided into several steps.

\ss

{\it Step 1.} First of all, any neighborhood $\cO$ of $G_0$ in $G$
can be covered by a finite number of the images of the following
open subset of $\dbR^n$
 $$
 G'=\Big\{(w_1,\cdots,w_n)\in
 \dbR^n\;\Big|\;0<w_n<1-\sum_{i=1}^{n-1}w_i^2\Big\}
 $$
under diffeomorphisms $x_j=x_j(w_1,\cdots,w_n)$ ($1\le j\le n$) of
class $C^2(\cl{G'})$ so that the image of $\pa G'\cap\{w_n=0\}$ is
contained in $\pa G_0$. Such diffeomorphisms change the
coefficients of the parabolic operator $\cF$ in (\ref{hh6.1}), but
do not change its parabolicity and the adoptedness of solutions.
Therefore, it suffices to consider $G=G'$. Note also that those
diffeomorphisms do not change the time variable. Hence, to
simplify the notations and noting that the original $z$ vanishes
in $(0,T)\t G_0\t\O$, we may assume the resulting parabolic
equation in $(0,T)\t G'\t\O$ reads
 $$
 \left\{
 \ba{ll}
 \ds\cF z\equiv dz-\sum_{i,j}(a^{ij}z_i)_jdt=[\lan a,\n z\ran+bz]dt+cz
 dw(t),\qq\hb{ in }(0,T)\t G'\t\O,\\
 \ns
 \supp z\subset(0,T)\t\big\{(x',x_n)\;\big|\;x_n\ge 0\big\}\t\O,
 \ea
 \right.
 $$
where $x'=(x_1,\cdots,x_{n-1})$ and $x=(x',x_n)$.

Next, we introduce a ``partial Holmgren coordinate transform"
$F:\; G'\to\dbR^n$ as follows:
 \bel{pa1}
 \left\{
 \ba{ll}
 \tilde x'=x',\\
 \tilde x_n=|x'|^2+x_n.
 \ea
 \right.
 \ee
In is easy to see that
 $$
 F(G')=\Big\{(\tilde x',\tilde x_n)\;\Big|\;|\tilde x'|^2<\tilde x_n<1\Big\}.
 $$
Again, the coordinate transform $F$ does not change the
parabolicity of $\cF$  and the adoptedness of solutions. Hence, to
simplify the notations, we may assume the resulting parabolic
equation to be the following:
 \bel{he4}
 \left\{
 \ba{ll}
 \ds\cF z\equiv dz-\sum_{i,j}(a^{ij}z_i)_jdt=[\lan a,\n z\ran+bz]dt+cz
 dw(t),\qq\hb{ in }(0,T)\t U\t\O,\\
 \ns
 \supp z\subset(0,T)\t\big\{(x',x_n)\;\big|\; x_n\ge |
 x'|^2\big\}\t\O,
 \ea
 \right.
 \ee
where $U=\big\{(x',x_n)\;\big|\;|x'|^2< x_n<1\big\}$. It suffices
to show that
 \bel{oohe13}
 z\equiv 0,\qq \hb{in }(0,T)\t U\t\O.
 \ee

 Finally, fix any $r_0$ and $r_1$ with $0<r_0<r_1<1$, we choose a function $\rho\in C^\i[0,1]$ so that
  \bel{he6}
  \left\{
  \ba{ll}
  0\le\rho(x_n)\le 1, \qq & x_n\in [0,1],\\
  \rho(x_n)\equiv 1,&0\le x_n\le r_0,\\
  \rho(x_n)\equiv 0,& r_1\le x_n\le 1.
  \ea
  \right.
  \ee
 Put
  \bel{he8}
  u=u(t,x',x_n)\=\rho(x_n)z(t,x',x_n),\qq (t,x)\in (0,T)\t U\t\O.
  \ee
 Then, by the first equation in (\ref{he4}), we have
   \bel{00he4}
   \3n\ba{ll}
\ds du-\sum_{i,j}(a^{ij}u_i)_jdt\\
\ns \ds=d(\rho
z)-\sum_{i,j}\big(a^{ij}(\rho z)_i\big)_jdt\\
 \ns
 \ds=\rho \cF z-\sum_{i,j}\[(a^{ij}\rho_iz)_j+a^{ij}\rho_jz_i\]dt\\
 \ns
 \ds=\Big\{\rho (\lan a,\n
z\ran+bz)-\sum_{i,j}\[(a^{ij}\rho_iz)_j+a^{ij}\rho_jz_i\]\Big\}dt+\rho
cz
 dw(t),\ \hb{ in }(0,T)\t U\t\O;
  \ea
 \ee
while, by the second equation in (\ref{he4}) and noting
 (\ref{he6}), one has
  \bel{he7}
  u=0,\qq \hb{on } (0,T)\t\pa U\t\O.
  \ee

 \ss

 {\it Step 2.} The above transforms do not change the adaptedness of $z$, and hence that of $u$. We now apply Theorem \ref{c1t3} to $u$ given by
 (\ref{he8}), $Q$ replaced by $(0,T)\t U$, and
  \bel{he12}
  \psi=\psi(x)=1-x_n,\qq x\in \cl U.
  \ee
 By (\ref{h7}) in Theorem \ref{c1t3}, and
 noting (\ref{00he4}), we conclude that there is a constant $C>0$ such that
 for any sufficiently large $\l$ and $\mu$, it holds
  \bel{he10}
 \ba{ll}\ds
 \l^3\mu^4\mathbb{E}\int_0^T\int_U\f^3\th^2u^2dxdt+\l\mu^2\mathbb{E}\int_0^T\int_U\f\th^2|\n
 u|^2dxdt\\
  \ns
  \ds\le C\[\mathbb{E}\int_0^T\int_U \th^2\Big\{\rho (\lan a,\n
z\ran+bz)-\sum_{i,j}\[(a^{ij}\rho_iz)_j+a^{ij}\rho_jz_i\]\Big\}^2dxdt\\
  \ns
  \ds\q+
 \mathbb{E}\int_0^T\int_U \th^2\sum_{i,j}
 a^{ij}(\rho cz)_i(\rho cz)_jdxdt\\
 \ns
 \ds\q+\mathbb{E}\int_0^T\int_U \th^2\[A-\sum_{i,j}
 \(b^{ij}\ell_i\ell_j+(b^{ij}\ell_i)_j\)\] (\rho cz)^2dxdt\],
 \ea
 \ee
where $\ds A=-\sum_{i,j}
 (a^{ij}\ell_i\ell_j-a^{ij}_j\ell_i
 +a^{ij}\ell_{ij})$.

 By the first estimate in (\ref{h6}) and noting our assumptions on $a$, $b$ and $c$,  we get
 \bel{he11}
  \3n\ba{ll}
  \ds\mathbb{E}\int_0^T\int_U \th^2\Big\{\rho (\lan a,\n
z\ran+bz)-\sum_{i,j}\[(a^{ij}\rho_iz)_j+a^{ij}\rho_jz_i\]\Big\}^2dxdt\\
  \ns
  \ds\q+
 \mathbb{E}\int_0^T\int_U \th^2\sum_{i,j}
 a^{ij}(\rho cz)_i(\rho cz)_jdxdt\\
 \ns
 \ds\q+\mathbb{E}\int_0^T\int_U \th^2\[A-\sum_{i,j}
 \(b^{ij}\ell_i\ell_j+(b^{ij}\ell_i)_j\)\](\rho cz)^2dxdt\\
   \ns
  \ds\le C\mathbb{E}\int_0^T\int_U\th^2(\l^2\mu^2\f^2z^2+|\n z|^2)dxdt.
  \ea
  \ee

 On the other hand, by (\ref{he6}) and (\ref{he8}), one finds
  \bel{oohe12}
  \ba{ll}
  \ds\l\mu^2\mathbb{E}\int_0^T\int_U\f\th^2|\n
  u|^2dxdt+\l^3\mu^4\mathbb{E}\int_0^T\int_U\f^3\th^2u^2dxdt\\
   \ns
  \ds\ge
  \mathbb{E}\int_0^T\int_{U\cap\{0<x_n<r_0\}}\th^2\[\l\mu^2\f|\n
  z|^2+\l^3\mu^4\f^3z^2\]dxdt.
  \ea
  \ee
Hence, combining (\ref{he10})--(\ref{oohe12}), and choosing $\l$
and $\mu$ large enough, we arrive at
 \bel{he13}
  \ba{ll}
  \ds\mathbb{E}\int_0^T\int_{U\cap\{0<x_n<r_0\}}\th^2\[\l\mu^2\f|\n
  z|^2+\l^3\mu^4\f^3z^2\]dxdt\\
   \ns
  \ds\le
  C\mathbb{E}\int_0^T\int_{U\cap\{r_0<x_n<1\}}\th^2(\l^2\mu^2\f^2z^2+|\n z|^2)dxdt.
  \ea
  \ee

  \ss

  {\it Step 3.} From now on, we fix $\mu$. Also, we fix any $\kappa_1\in (0,r_0)$. Noting the definition of $\psi$  in (\ref{he12}) implies
that $\th=\th(t,x_n)$ is decreasing with respect to $x_n$, from
(\ref{he13}), we deduce that
 \bel{pphe14}
  \ba{ll}
  \ds\l^3\mu^4\mathbb{E}\int_0^T\int_{U\cap\{0<x_n<\kappa_1\}}|\th(t,\kappa_1)|^2\f^3z^2dxdt\\
   \ns
  \ds\le \mathbb{E}\int_0^T\int_{U\cap\{0<x_n<\kappa_1\}}|\th(t,\kappa_1)|^2\[\l\mu^2\f|\n
  z|^2+\l^3\mu^4\f^3z^2\]dxdt\\
   \ns
  \ds\le \mathbb{E}\int_0^T\int_{U\cap\{0<x_n<\kappa_1\}}\th^2\[\l\mu^2\f|\n
  z|^2+\l^3\mu^4\f^3z^2\]dxdt\\
   \ns
  \ds\le \mathbb{E}\int_0^T\int_{U\cap\{0<x_n<r_0\}}\th^2\[\l\mu^2\f|\n
  z|^2+\l^3\mu^4\f^3z^2\]dxdt\\
   \ns
  \ds\le C\mathbb{E}\int_0^T\int_{U\cap\{r_0<x_n<1\}}\th^2(\l^2\mu^2\f^2z^2+|\n z|^2)dxdt\\
   \ns
  \ds\le
  C\mathbb{E}\int_0^T\int_{U\cap\{r_0<x_n<1\}}|\th(t,r_0)|^2(\l^2\mu^2\f^2z^2+|\n z|^2)dxdt,
  \ea
  \ee
 for a constant $C>0$, independent of $\l$.

Further, fix any $\kappa_2\in (0,1)$. Noting  that
$\th=\th(t,x_n)$ is increasing (\resp decreasing) with respect to
$t$ in $[0,T/2]$ (\resp $(T/2,1]$), we deduce that
 \bel{ann1}
 \ba{ll}
 \ds\mathbb{E}\int_0^T\int_{U\cap\{0<x_n<\kappa_1\}}|\th(t,\kappa_1)|^2\f^3z^2dxdt\\\ns
 \ds\ge
 \mathbb{E}\int_{(1-\kappa_2)T/2}^{(1+\kappa_2)T/2}\int_{U\cap\{0<x_n<\kappa_1\}}|\th(t,\kappa_1)|^2\f^3z^2dxdt\\\ns
 \ds\ge|\th((1-\kappa_2)T/2,\kappa_1)|^2
 \mathbb{E}\int_{(1-\kappa_2)T/2}^{(1+\kappa_2)T/2}\int_{U\cap\{0<x_n<\kappa_1\}}\f^3z^2dxdt,
 \ea
 \ee
and
 \bel{ann2}
 \ba{ll}
 \ds
 \mathbb{E}\int_0^T\int_{U\cap\{r_0<x_n<1\}}|\th(t,r_0)|^2(\l^2\mu^2\f^2z^2+|\n
 z|^2)dxdt\\\ns\ds\le |\th(T/2,r_0)|^2\mathbb{E}\int_0^T\int_{U\cap\{r_0<x_n<1\}}(\l^2\mu^2\f^2z^2+|\n
 z|^2)dxdt.
 \ea
 \ee
Combining (\ref{pphe14})--(\ref{ann2}), we end up with
 \bel{ann7}
 \ba{ll}
 \ds\l^3\mu^4|\th((1-\kappa_2)T/2,\kappa_1)|^2
 \mathbb{E}\int_{(1-\kappa_2)T/2}^{(1+\kappa_2)T/2}\int_{U\cap\{0<x_n<\kappa_1\}}\f^3z^2dxdt\\\ns\ds\le C
 |\th(T/2,r_0)|^2\mathbb{E}\int_0^T\int_{U\cap\{r_0<x_n<1\}}(\l^2\mu^2\f^2z^2+|\n
 z|^2)dxdt.
 \ea
 \ee

By (\ref{h5}), (\ref{alphad}) and (\ref{he12}), we find
 \bel{ann3}\ba{ll}
 \ds
 |\th((1-\kappa_2)T/2,\kappa_1)|^2=\exp\left\{{2\l[e^{(1-\kappa_1)\mu}-e^{2\mu}]4^k\over
 (1-\kappa_2^2)^kT^{2k}}\right\},\\\ns\ds |\th(T/2,r_0)|^2=\exp\left\{{2\l[e^{(1-r_0)\mu}-e^{2\mu}]4^k\over
 T^{2k}}\right\}.
 \ea
 \ee
We now choose $\kappa_2$ to be
 \bel{ann4}
 \kappa_2=\sqrt{1-\sqrt[k]{{e^{2\mu}-e^{(1-\kappa_1)\mu}\over
 e^{2\mu}-e^{(1-r_0)\mu}}}}.
 \ee
Since $\kappa_1\in (0,r_0)$, one sees that $\kappa_2\in (0,1)$.
Moreover, by (\ref{ann4}), we have
 \bel{ann5}
 {e^{(1-\kappa_1)\mu}-e^{2\mu}\over
 (1-\kappa_2^2)^k}=e^{(1-r_0)\mu}-e^{2\mu}.
 \ee
Combining (\ref{ann3}) and (\ref{ann5}), it follows that
  \bel{ann6}
 |\th((1-\kappa_2)T/2,\kappa_1)|^2= |\th(T/2,r_0)|^2.
 \ee

 Now, by (\ref{ann7}) and noting (\ref{ann6}), we conclude that
 \bel{ann8}
 \ba{ll}
 \ds\l^3\mu^4
 \mathbb{E}\int_{(1-\kappa_2)T/2}^{(1+\kappa_2)T/2}\int_{U\cap\{0<x_n<\kappa_1\}}\f^3z^2dxdt\\\ns\ds\le C
 \mathbb{E}\int_0^T\int_{U\cap\{r_0<x_n<1\}}(\l^2\mu^2\f^2z^2+|\n
 z|^2)dxdt.
 \ea
 \ee

 Letting
 $\l\to+\i$ in (\ref{ann8}), we conclude that
  $$
  z\equiv 0, \qq \hb{in }
 \big((1-\kappa_2)T/2,(1+\kappa_2)T/2\big)\t(U\cap\{0<x_n<\kappa_1\})\t \O.
 $$
Hence,
 $$
  z(T/2,\cd)\equiv 0, \qq \hb{in }(U\cap\{0<x_n<\kappa_1\})\t \O.
 $$
Since $r_0$ (\resp $\kappa_1$) can be chosen as close to $1$
(\resp $r_0$) as one likes, one concludes that
  $$
  z(T/2,\cd)\equiv 0, \qq \hb{in }U\t \O.
  $$
Replace $T$ by any given $t_0\in (0,T)$. Then, the above argument
yields $z(t_0/2,\cd)\equiv 0$ in $U\t \O$. Hence,
 $$
  z\equiv 0, \qq \hb{in }(0,T/2]\t U\t \O.
  $$
 Applying this argument to $z(\cd+T/2,\cd)$, it follows that
 $$
  z\equiv 0, \qq \hb{in }(T/2,3T/4]\t U\t \O.
  $$
Repeating this procedure, we arrive at (\ref{oohe13}). This
completes the proof of Theorem \ref{1t3}.
\endpf

\end{document}